\documentclass[journal]{IEEEtran}
\usepackage{amsmath}
\usepackage[pdftex]{graphicx}
\usepackage{amssymb}
\usepackage{booktabs}
\usepackage{cite}
\usepackage{color}
\usepackage{enumerate}
\usepackage{subcaption}

\newtheorem{theorem}{Theorem}

\def\tinyL{{\tiny \mbox{L}}}
\def\tinyU{{\tiny \mbox{U}}}

\allowdisplaybreaks

\begin{document}

%\title{Affinely Adjustable Robust Power Dispatch Integrating Do-Not-Exceed Limits of Wind Power}
\title{Distributionally Robust Co-Optimization of Power Dispatch and Do-Not-Exceed Limits}
%\author{,~\IEEEmembership{Member,~IEEE,}
%        ,~\IEEEmembership{Fellow,~OSA,}
%      ,~\IEEEmembership{Life~Fellow,~IEEE}%

\author{Hongyan~Ma,~\IEEEmembership{}
        Ruiwei~Jiang,~\IEEEmembership{Member,~IEEE,}
       and~Zheng~Yan~\IEEEmembership{}
%\thanks{Manuscript received April 19, 2018; revised August 26, 2018.}
\thanks{H. Ma and Z. Yan are with the Key Laboratory of Control of Power Transmission and Conversion, Ministry of Education, Department of Electrical Engineering, Shanghai Jiao Tong University, Shanghai 200240, China (e-mail: hahaha91644@sjtu.edu.cn).}
\thanks{R. Jiang is with the Department of Industrial and Operations Engineering, University of Michigan, Ann Arbor, MI 48109, USA (e-mail: ruiwei@umich.edu).}
\thanks{This work is supported in part by the U.S. National Science Foundation (CMMI-1662774) and the National Key R\&D Program of China (Technology and application of wind power/photovoltaic power prediction for promoting renewable energy consumption, 2018YFB0904200).}
}

%\markboth{Journal of Class Files,~Vol.~14, No.~8, August~2015}%

\maketitle

\begin{abstract}
To address the challenge of the renewable energy uncertainty, the ISO New England (ISO-NE) has proposed to apply do-not-exceed (DNE) limits, which represent the maximum nodal injection of renewable energy the grid can accommodate. Unfortunately, it appears challenging to compute DNE limits that simultaneously maintain the system flexibility and incorporate a large portion of the available renewable energy at the minimum cost. In addition, it is often challenging to accurately estimate the joint probability distribution of the renewable energy. In this paper, we propose a two-stage distributionally robust optimization model that co-optimizes the power dispatch and the DNE limits, by adopting an affinely adjustable power re-dispatch and an adjustable joint chance constraint that measures the renewable utilization. Notably, this model admits a second-order conic reformulation that can be efficiently solved by the commercial solvers (e.g., MOSEK). We conduct case studies based on modified IEEE test instances to demonstrate the effectiveness of the proposed approach and analyze the trade-off among the system flexibility, the renewable utilization, and the dispatch cost.
%Due to the challenge of uncertainty from renewable resources, New England ISO proposes the do-not-exceed limit which is the maximum renewable generation ranges that the system can accommodate. In this paper, we propose a robust power dispatch approach through incorporating an affinely adjustable corrective re-dispatch and integrating the do-not-exceed limit of wind power, which co-optimizes the power pre-dispatch strategies and the admissible ranges of wind power outputs. The probability that the wind power can be covered by the admissible ranges is adopted to measure the utilization of wind power in the objective. This approach is formulated as a two-stage robust optimization with the uncertainty set as decision variable, and we propose the algorithms of approximating the adjustable uncertainty set and joint chance constraints to reduce the computational scale. Case studies on the modified IEEE systems display the frontier of the admissible ranges, confidence level and the dispatch cost and compare the performance of the proposed approach with the original do-not-exceed limit approach. It turns out that the proposed formulation and algorithm is effective and the proposed approach is robust and able to utilize more wind power and reduce operation costs.
\end{abstract}


\IEEEpeerreviewmaketitle

\section*{Nomenclature}
\addcontentsline{toc}{section}{Nomenclature}
\begin{IEEEdescription}[\IEEEusemathlabelsep\IEEEsetlabelwidth{$V_1,V_2$}]
	\item [\textcolor{black}{\bf Indices and Sets}]
	\item[$t,i,k,n,l$] \hspace{0.2cm} Index for time period, thermal unit, renewable resource, node, and transmission line, respectively.
    \item[$T,I,K,N,L$] \hspace{0.6cm} Numbers of time periods, thermal units, renewable resources, nodes, and transmission lines, respectively.
%	\item[$n(i), n(k)$] \hspace{0.2cm} Index for the node where thermal unit $i$ and renewable resource $k$ is located, respectively. ({\color{red}may not be needed after using abstract notation})
	\item[{$[i(n)], [k(n)]$}] \hspace{0.6cm} Sets of thermal units and renewable resources at node $n$, respectively.
	\item[{$[M]$}] $[M]:=\{1, \ldots, M\}$ for positive integer $M$.
	\vspace{0.2cm}
	
	\item [\textcolor{black}{\bf Parameters}]
    \item[$C_i(\cdot)$] Fuel cost function of thermal unit $i$.
    \item[$c^+_{kt}, c^-_{kt}$] Unit cost of overestimating and underestimating the output of renewable resource $k$ during time period $t$, respectively.
    \item[$\hat{w}_{kt}$]   Forecasted output of renewable resource $k$ during time period $t$.
    \item[$d_{nt}$]     Load of node $n$ during time period $t$.
    \item[$\bar{F_l}$]  Transmission capacity limit of line $l$.
    \item[$f_{nl}$]     Shift distribution factor of node $n$ to line $l$.
    \item[$p_i^{\min},p_i^{\max}$]  \hspace{0.3cm} Minimum and Maximum generation capacity of thermal unit $i$, respectively.
    \item[$r_i^\text{up},r_i^\text{dn}$]  Upward and downward ramp-rate of thermal unit $i$ (MW/min), respectively.
    \item[$\Delta_\text{d},\Delta_\text{r}$]  Dispatch interval (min) and response time window (min), respectively. %({\color{red}Why subscript $t$? unit = min or hr? in the case study, it says 1hr and 5 min, resp.})
    \item[$w_k^{\min},w_k^{\max}$]  \hspace{0.5cm} Minimum and maximum generation capacity of renewable resource $k$, respectively.
    \item[$\mu_{kt},\sigma_{kt}$]  Empirical mean and empirical variance of the prediction error of renewable output $k$ during time period $t$, respectively.
    \item[$\delta$]  Weight coefficient between power dispatch and renewable utilization.
    \item[$\varepsilon_{kt}$]  Random variable representing the output deviation of renewable resource $k$ from its forecast value $\hat{w}_{kt}$ during time period $t$.
    \item[$u_0$] Lower bound on renewable utilization.
    \vspace{0.2cm}
	
	\item [\textcolor{black}{\bf Decision Variables}]
    \item[$\hat{p}_{it}$]  Scheduled generation amount of thermal unit $i$ during time period $t$.
    \item[$p_{it}$]  Actual generation amount of thermal unit $i$ at time $t$.
    \item[$\varepsilon_{kt}^{\tinyL},\varepsilon_{kt}^{\tinyU}$]    Lower and upper do-not-exceed limits of $\varepsilon_{kt}$, respectively.
    \item[$u$] Renewable utilization probability.
    \item[$B_{ikt},b_{ikt}$] \hspace{0.2cm} Coefficients of the affine decision rule.
%    \item[$Q_{kt},q_{kt}$] \hspace{0.2cm} Auxiliary variables.
%    \item[$H_{ikt},h_{ikt}$] \hspace{0.2cm} Auxiliary variables.
%    \item[$\lambda_{klt}^\text{F},\gamma_{klt}^\text{F}$] \hspace{0.2cm} Auxiliary variables.
%    \item[$\lambda_{ikt}^\text{R},\gamma_{ikt}^\text{R}$] \hspace{0.2cm} Auxiliary variables.
%    \item[$\lambda_{ikt}^\text{C},\gamma_{ikt}^\text{C}$] \hspace{0.2cm} Auxiliary variables.
    \item[$r_{kt},s_{kt},z_{kt}$] \hspace{0.5cm} Auxiliary dual variables in the reformulation of the adjustable joint chance constraint.
\end{IEEEdescription}

\section{Introduction}

\IEEEPARstart{T}{he} renewable energy (e.g., wind and solar power) leads to random nodal injections in the power grid and presents a significant challenge to the power system operation. Many methods have been proposed to hedge against the renewable energy uncertainty, including stochastic programming (see, e.g.,~\cite{wu2007stochastic,papavasiliou2011stochastic,wang2012stochastic}) and robust optimization approaches (see, e.g.,~\cite{jiang2012robust,bertsimas2013adaptive}). These approaches incorporate the uncertainty based on pre-specified models, e.g., probability distribution in stochastic programming (SP) approaches and uncertainty set in robust optimization (RO) approaches. Although SP and RO approaches are widely applied, we may still need to address the following two challenges:
\begin{description}
\item[{\bf Challenge 1}] \hspace{1cm}It may be too costly or even infeasible to treat the renewable energy as non-dispatchable resource and balance its variation by regulating other dispatchable resources, especially when the renewable penetration is high~\cite{zhao2015dnelimit,wei2015dispatchable,shao2017security}.
\item[{\bf Challenge 2}] \hspace{1cm}It is often challenging to accurately estimate the joint probability distribution of the renewable energy. Consequently, the solution obtained from a SP model can perform worse in out-of-sample tests than in the in-sample tests (see, e.g.,~\cite{zymler2011distributionally,wiesemann2014distributionally,jiang2016data}).
\end{description}
%On one hand, there would has no feasible solution under the predefined uncertainty set if the system flexibility is not enough or the wind power deviation is too large. On the other hand, even if the flexibility is sufficient, the robust solution might be very expensive because the flexible resources might be employed to reserve ramp capabilities to deal with the uncertainty at high electricity prices.

To address Challenge 1, the ISO-NE proposes an inspiring concept of do-not-exceed (DNE) limits under a given power dispatch strategy~\cite{zhao2015dnelimit}. The DNE limits assign an admissible range of renewable energy to each node of the transmission system. This provides a clear guideline for utilizing renewable energy: the system accommodates any nodal injection that lies within the admissible range, and otherwise emergency regulations (e.g., renewable energy curtailment, fast-starting units, and load shedding) may have to be used. In addition, the DNE limits also offer a convenient way of defining and measuring of the system flexibility~\cite{zhao2016framework}. %Recently, the DNE limits have received increasing attention in the literature. \cite{wei2015dispatchable,wei2015real} show that the admissible range of a power grid is mathematically equivalent to a polytope. This generalizes the DNE limits, which take the form of a hypercube.

Recently, the DNE limits have received increasing attention in the literature. \cite{wei2015dispatchable,wei2015real} show that the admissible range of a power grid is mathematically equivalent to a polytope. This generalizes the concept of the DNE limits, which take the form of a hypercube. In~\cite{zhao2015dnelimit,wei2015dispatchable,wei2015real}, the admissible range is obtained based on a given power dispatch strategy, which, however, may not be optimal for accommodating renewable generation. As a result, this might underestimate the dispatch capability of the power system in accommodating renewable energy. As an alternative, many studies propose to co-optimize the power dispatch and the DNE limits. \cite{li2015robust} proposes a single-stage RO model that co-optimizes the power dispatch and the polytopic admissible range. To incorporate recourse actions (e.g., power re-dispatch),~\cite{li2015adjustable} proposes an adjustable RO model in which recourse actions follow an affine decision rule (ADR) with given coefficients. \cite{wei2016dispatchability} studies an adjustable RO model with ADR and optimized coefficients. Additionally, the proposed model in~\cite{wei2016dispatchability} incorporates risk criteria based on the radius and the coverage probability of the admissible range. Without applying an ADR,~\cite{shao2017power} considers an adjustable RO model that co-optimizes power dispatch and DNE limits with full recourse. Later,~\cite{shao2017security} extends~\cite{shao2017power} by incorporating unit commitment (UC) into the co-optimization, and employs the column-and-constraint generation (CCG) algorithm~\cite{zeng2013solving} to solve the proposed model. \cite{wang2017robustrcuc} also considers an adjustable RO model with full recourse that incorporates UC and a polytopic admissible range. Furthermore,~\cite{wang2016risk} considers the risk of the renewable energy being realized outside of the admissible region, which result in, e.g., curtailment of renewable energy. Differently,~\cite{qiu2017data} models this risk by maximizing the probability that the renewable energy being realized within the DNE limits. Then, this model is solved by using the sample average approximation algorithm. \cite{li2016multi} considers an adjustable RO model with both discrete and continuous recourse and proposes to solve this model with a nested CCG algorithm. \cite{korad2015zonal} incorporates topology control into the co-optimization and considers zonal DNE limits. It is worth mentioning that solving the adjustable RO model with full recourse requires repeatedly solving mixed-integer programs with big-M coefficients (see, e.g.,~\cite{wei2015real,korad2015zonal,wang2016risk,shao2017security}), which may be challenging when many nodes of the power grid incorporate renewable energy.

A natural way of mitigating Challenge 2 is to employ distributionally robust optimization (DRO). In contrast to SP that considers a single probability distribution, DRO considers a family of probability distributions that are plausible of modeling the renewable energy. We term the family of distributions as an ambiguity set. In the existing literature, ambiguity sets based on the moments of uncertainty (e.g., mean, variance, etc.) are commonly applied (see, e.g.,~\cite{xiong2017distributionally,zhang2017distributionally,xie2018distributionally,zhao2018distributionally}). Other distributional information based on, e.g., the Wasserstein distance~\cite{esfahani2017data,wang2017risk}, the $\phi$-divergence~\cite{jiang2016data,chen2018distributionally}, and the unimodality~\cite{wei2016dispatchability,li2016distributionally,li2017ambiguous}, have also been proposed to characterize the ambiguity set. Accordingly, DRO formulates a robust counterpart of SP and hedges against the worst-case probability distribution within the ambiguity set.

In this paper, we consider a distributionally robust (DR) co-optimization model for the power dispatch and the DNE limits. Our model follows~\cite{wang2016risk,qiu2017data} to incorporate the operational risks. Specifically, we consider the DR probability that the renewable energy being realized within the DNE limits. We further extend the model to incorporate the DR expected cost of overestimating/underestimating the renewable energy. By using ADR with optimized coefficients, we show that this model admits a conic programming reformulation that can be efficiently solved by the commercial solvers (e.g., MOSEK). The proposed model assumes a fixed UC and power grid topology. Nevertheless, these decisions can also be incorporated into this model with slight changes, leading to mixed 0-1 conic programming reformulations.

The remainder of this paper is organized as follows. Section \ref{sec:formulation} presents the mathematical formulation and Section \ref{sec:solution} describes the solution methodology. Section \ref{sec:dr-expected-cost} extends the model and solution methodology to incorporate alternative operational risks. Section \ref{sec:computation} reports the case studies that demonstrate the effectiveness of the proposed approach, before we draw conclusions in Section \ref{sec:conclusion}.

\section{Mathematical Formulation} \label{sec:formulation}
We describe the co-optimization model of power dispatch and DNE limits in Section \ref{sec:dne} and the adjustable DR chance constraint in Sections \ref{sec:risk}.

\subsection{DNE Limits} \label{sec:dne}
Given the forecasted renewable energy outputs $\hat{w}_{kt}$, the nominal economic dispatch (ED) model maintains the generation-load balance under operational restrictions. Mathematically, we formulate the constraints of a multi-period nominal ED model as follows:
\begin{subequations}
\begin{align}
%&\min_{\hat{p}_{it}} \sum_t\sum_iC_i(\hat{p}_{it}) \label{EDob}\\
%\text{s.t.}&\sum_i\hat{p}_{it}+\sum_k\hat{w}_{kt}=\sum_nd_{nt},\forall t \label{EDc1}\\
&\sum_{i \in [I]} \hat{p}_{it}+\sum_{k \in [K]} \hat{w}_{kt} = \sum_{n \in [N]} d_{nt}, \ \forall t \in [T], \label{EDc1}\\
& -\bar{F_l}\leq \sum_{n \in [N]} f_{nl} \left( \sum_{i \in [i(n)]} \hat{p}_{it} +\sum_{k \in [k(n)]} \hat{w}_{kt} - d_{nt} \right) \leq \bar{F_l}, \nonumber \\
& \forall l \in [L], \ \forall t \in [T], \label{EDc2} \\
& p_i^\text{min} \leq \hat{p}_{it} \leq p_i^\text{max}, \ \forall i \in [I], \ \forall t \in [T], \label{EDc3} \\
& -r_i^\text{dn}  \Delta_\text{d} \leq \hat{p}_{it}-\hat{p}_{i,t-1}\leq r_i^\text{up}  \Delta_\text{d}, \ \forall i \in [I], \ \forall t \in [T], \label{EDc4}
\end{align}
\end{subequations}
where $\hat{p}_{it}$ represents the pre-dispatch strategy based on the forecast renewable generation $\hat{w}_{kt}$, constraints \eqref{EDc1} represent the generation-load balance, \eqref{EDc2} represent the transmission line capacity restrictions based on the dc approximation of the power flow equations, \eqref{EDc3} represent the capacity limits of the thermal units, and \eqref{EDc4} represent the ramp-rate limits of the thermal units. When taking the uncertainty of renewable energy and the DNE limits into account, the power system aims to accommodate any nodal injections of renewable energy through corrective power re-dispatch, as long as such injections lie within the DNE limits. We formulate this requirement as follows for all $t \in [T]$:
\begin{subequations}
\begin{align}
& \forall \varepsilon_t \in [\varepsilon_t^{\tinyL},\varepsilon_t^{\tinyU}], \ \mbox{there exist } \{p_{it}(\varepsilon_t)\}_{i \in [I]} \mbox{ such that:} \notag \\
& \sum_{i \in [I]} p_{it}(\varepsilon_t) + \sum_{k \in [K]}(\hat{w}_{kt}+\varepsilon_{kt}) = \sum_{n \in [N]}d_{nt}, \label{MOc4}\\
& -\bar{F_l}\leq \sum_{n \in [N]} f_{nl} \left( \sum_{i \in [i(n)]} p_{it}(\varepsilon_t)+ \sum_{k \in [k(n)]} (\hat{w}_{kt} +\varepsilon_{kt})\right) \notag \\
&  \leq \bar{F_l}, \ \forall l \in [L], \label{MOc5} \\
& p_i^\text{min} \leq  p_{it}(\varepsilon_t) \leq p_i^\text{max}, \ \forall i \in [I], \label{MOc6}\\
& -r_i^\text{dn} \Delta_\text{d} \leq p_{it}(\varepsilon_t)-p_{i,t-1}(\varepsilon_{t-1})\leq r_i^\text{up} \Delta_\text{d}, \ \forall i \in [I], \label{MOc7} \\
& -r_i^\text{dn} \Delta_\text{r} \leq p_{it}(\varepsilon_t)-\hat{p}_{it} \leq r_i^\text{up} \Delta_\text{r}, \ \forall i \in [I], \label{MOc8}
\end{align}
where constraints \eqref{MOc4}--\eqref{MOc7} are counterparts of \eqref{EDc1}--\eqref{EDc4} with regard to the power re-dispatch variables $p_{it}(\varepsilon_t)$ and constraints \eqref{MOc8} represent the ramping capacity limits within the response time window. In this paper, we assume that $p_{it}(\varepsilon_t)$ follows an ADR, i.e., $p_{it}(\varepsilon_t)$ is the following affine function of $\varepsilon_t$:
\begin{equation}
p_{it}(\varepsilon_{t}) = \hat{p}_{it}+\sum_k(B_{ikt}\varepsilon_{kt}+b_{ikt}), \ \forall i \in [I], \ \forall t \in [T], \label{MOc9}
\end{equation}
\end{subequations}
where $B_{ikt}$ and $b_{ikt}$ represent the response of $p_{it}(\varepsilon_{t})$ to the forecast error $\varepsilon_t$ and can be adjusted to optimize the objective function (to be specified in Section \ref{sec:risk}). On the one hand, the ADR corresponds to the incremental output of the automatic generation control (AGC) units as an affine function of the renewable generation deviation. For the non-AGC units, we can set $B_{ikt} = b_{ikt} = 0$. On the other hand, the ADR restricts the search space of the recourse variables $p_{it}(\varepsilon_{t})$ and so yields a conservative approximation of the constraints \eqref{MOc4}--\eqref{MOc8}.

\subsection{Adjustable DR Joint Chance Constraints} \label{sec:risk}
We note that formulations \eqref{EDc1}--\eqref{EDc4} and \eqref{MOc4}--\eqref{MOc8} do not incorporate any distributional information of the forecast error $\varepsilon_{kt}$. This may cause a mismatch between the DNE limits and renewable energy. For example, it may be unlikely that the renewable generation is realized within the DNE limits, and accordingly we may curtail a significant portion of the renewable generation. To address this challenge, we first designate that the DNE limits contain the forecasted output of renewable energy and lie within the capacity limits of the renewable generation:
\begin{align}
& w^{\min}_k \leq \hat{w}_{kt}+\varepsilon_{kt}^{\tinyL} \leq \hat{w}_{kt}\leq \hat{w}_{kt}+\varepsilon_{kt}^{\tinyU} \leq w^{\max}_k, \nonumber \\
& \forall k \in [K], \ \forall t \in [T]. \label{MOc10}
\end{align}
Second, we consider an \emph{adjustable} joint chance constraint to measure the utilization of renewable energy:
\begin{subequations}
\begin{align}
%&\max_{\varepsilon_{kt}^\text{U},\varepsilon_{kt}^\text{L}}  a \label{Pob}\\
%\text{s.t.} &\mathbb{P}(\varepsilon _t\in [\varepsilon_t^{\tinyL},\varepsilon_t^{\tinyU}]) \geq a,\forall t \label{Pc2}\\
& \inf_{\mathbb{P} \in \mathcal{D}}\mathbb{P} \left(\varepsilon_t \in [\varepsilon_t^{\tinyL},\varepsilon_t^{\tinyU}] \right) \geq u, \ \forall t \in [T], \label{Pc2}\\
& u_0 \leq u \leq 1, \label{Pc3}
\end{align}
\end{subequations}
where $u$ estimates the probability of fully utilizing the renewable energy and $u_0$ represents a lower bound of $u$. In this paper, we assume that $u_0 > 2/3$ (see Theorem~\ref{theorem-1}). This assumption is not very restrictive because power system operators often desire high utilization of renewable energy (see, e.g.,~\cite{fink2009wind}). Additionally, we note that $u$ represents a decision variable in our model and can be adjusted to optimize the trade-off between the power dispatch cost and the renewable utilization. In addition, we consider an ambiguity set $\mathcal{D}$ consisting of probability distributions $\mathbb{P}$ that (i) match the empirical mean $\mu_{kt}$ and empirical variance $\sigma_{kt}$ of each $\varepsilon_{kt}$ and (ii) is unimodal about $\mu_{kt}$, i.e.,
\begin{equation}
\mathcal{D} := \left\{\mathbb{P}: \hspace{-0.2cm} \begin{array}{l} \mathbb{E}_{\mathbb{P}}[\varepsilon_{kt}]=\mu_{kt}, \ \text{Var}(\varepsilon_{kt})=\sigma_{kt}^2, \\[0.2cm]
\varepsilon_{kt} \text{ is unimodal about } \mu_{kt}, \ \forall k \in [K], \forall t \in [T]
\end{array}\right\}. \label{ambiguity-set}
\end{equation}
Unimodality about $\mu_{kt}$ indicates that the probability density function of $\varepsilon_{kt}$, if exists, is nondecreasing from $0$ to $\mu_{kt}$ and is nonincreasing afterwards. In the literature, many probability distributions proposed for modeling the renewable energy forecast error are unimodal (see, e.g.,~\cite{doherty2005new,wang2008security,hodge2012wind}). It is worth mentioning that~\cite{wei2016dispatchability} consider the unimodality of the \emph{joint} probability distribution $\mathbb{P}$ of all $\varepsilon_{kt}$. In contrast, the unimodality in $\mathcal{D}$ is with respect to the \emph{marginal} distribution of each $\varepsilon_{kt}$, which is weaker than the joint unimodality assumed in~\cite{wei2016dispatchability} and easier to verify by the historical data. In addition, the ambiguity set $\mathcal{D}$ leads to a polynomially solvable reformulation (see Section \ref{sec:solution}).

We close this section by formulating the DR co-optimization (DRCO) model of power dispatch and DNE limits:
\begin{subequations}
\begin{align}
\min_{\substack{\hat{p},B,b,\\ \varepsilon^{\tinyL},\varepsilon^{\tinyU},u}} & \ \sum_{t \in [T]} \sum_{i \in [I]} C_i(\hat{p}_{it}) - \delta u \label{drco-obj} \\
\text{s.t.} \ & \ \mbox{\eqref{EDc1}--\eqref{EDc4}}, \ \mbox{\eqref{MOc4}--\eqref{MOc9}}, \ \eqref{MOc10}, \ \mbox{\eqref{Pc2}--\eqref{Pc3}}, \label{drco-con}
\end{align}
\end{subequations}
where $\delta$ represents the weight on the renewable utilization $u$. The system operator can set $\delta$ based on her trade-off between the renewable utilization and the dispatch cost. If $\delta$ is close to zero then the dispatch cost and the renewable utilization are low. As $\delta$ increases, both dispatch cost and renewable utilization increase. By gradually increasing the value of $\delta$ and re-solving model \eqref{drco-obj}--\eqref{drco-con}, we obtain a cost-utilization frontier that can clearly indicate the trade-off between these two performance measures (see Section \ref{sec:computation} for related case studies).

\section{Solution Methodology} \label{sec:solution}
We recast the DRCO model \eqref{drco-obj}--\eqref{drco-con} as a second-order conic program that is polynomially solvable. For notation brevity, we derive based on abstract notation. First, we represent constraints \eqref{MOc4}--\eqref{MOc9} in the following abstract form:
\begin{subequations}
\begin{align}
& \exists \ \varepsilon^{\tinyL}, \varepsilon^{\tinyU}, p(\varepsilon): \nonumber \\
& Tx + Wp(\varepsilon) \leq H\varepsilon, \ \forall \varepsilon \in [\varepsilon^{\tinyL}, \varepsilon^{\tinyU}], \label{ref-note-1} \\
& p(\varepsilon) = B\varepsilon + b, \label{ref-note-2}
\end{align}
where $p(\varepsilon)$ denotes the power re-dispatch variables, matrices $T$, $W$, and $H$ denote the given parameters in constraints \eqref{MOc4}--\eqref{MOc8}, and matrix $B$ and vector $b$ denote the variables in the ADR \eqref{MOc9}. Letting $E := \mbox{diag}(\varepsilon^{\tinyU} - \varepsilon^{\tinyL})$, we represent the hypercube $[\varepsilon^{\tinyL}, \varepsilon^{\tinyU}]$ as $\{\varepsilon^{\tinyL} + Ev: v \in [0, e]\}$, where $e$ denotes the vector of all ones.
%We note that this representation holds valid if $\varepsilon^{\tinyU}_{kt} - \varepsilon^{\tinyL}_{kt} > 0$ for all $k \in [K]$ and $t \ in [T]$, which holds valid as long as $u > 0$ at the optimality of model \eqref{drco-obj}--\eqref{drco-con}. Indeed, suppose that $\varepsilon^{\tinyU}_{kt} = \varepsilon^{\tinyL}_{kt}$ for some $k \in [K]$ and $t \in [T]$. Then, $\mathbb{P}(\varepsilon \in [\varepsilon^{\tinyL}, \varepsilon^{\tinyU}]) = 0$ for any $\mathbb{P}$ that is continuous (with respect to the Lebesgue measure). It follows that $\inf_{\mathbb{P} \in \mathcal{D}}\mathbb{P}(\varepsilon \in [\varepsilon^{\tinyL}, \varepsilon^{\tinyU}]) = 0$, which contradicts the assumption that $u > 0$.
Then, we recast \eqref{ref-note-1}--\eqref{ref-note-2} as
\begin{align}
& \exists \ \varepsilon^{\tinyL}, \varepsilon^{\tinyU}, B, b: \nonumber \\
& Tx + W(BEv + B\varepsilon^{\tinyL} + b\bigr) \leq H\varepsilon^{\tinyL} + HEv, \nonumber \\
& \forall v \in [0, e]. \label{ref-note-3}
\end{align}
We claim that \eqref{ref-note-3} is equivalent to
\begin{align}
& \exists \ \varepsilon^{\tinyL}, \varepsilon^{\tinyU}, S, s_0: \nonumber \\
& Tx + W(Sv + s_0) \leq H\varepsilon^{\tinyL} + HEv, \ \forall v \in [0, e]. \label{ref-note-4}
\end{align}
We now prove the equivalence $\mbox{\eqref{ref-note-3}} \Leftrightarrow \mbox{\eqref{ref-note-4}}$. One the one hand, suppose that there exist $B$ and $b$ such that \eqref{ref-note-3} holds valid. Then, we let $\bar{S} = BE$ and $\bar{s}_0 = B\varepsilon^{\tinyL} + b$ to yield $Tx + W(\bar{S}v + \bar{s}_0) \leq H\varepsilon^{\tinyL} + HEv$ for all $v \in [0, e]$. Hence, \eqref{ref-note-3} implies \eqref{ref-note-4}. On the other hand, suppose that there exist $S$ and $s_0$ such that \eqref{ref-note-4} holds valid. Then, we let $\bar{B} = SE^{-1}$ and $\bar{b}_0 = s_0 - SE^{-1}\varepsilon^{\tinyL}$ to yield $Tx + W(\bar{B}Ev + \bar{B}\varepsilon^{\tinyL} + \bar{b}_0) \leq H\varepsilon^{\tinyL} + HEv$ for all $v \in [0, e]$. Hence, \eqref{ref-note-4} also implies \eqref{ref-note-3}. Furthermore, we note that constraint \eqref{ref-note-4} holds valid if and only if $\sup_{v \in [0, e]} \{(WS - HE)v\} \leq H\varepsilon^{\tinyL} - Tx - Ws_0$, where the supremum operator is applied on each component of $(WS - HE)v$. Using the standard technique in robust optimization (see, e.g.,~\cite{soyster1973convex}), we recast this constraint, and so constraints \eqref{MOc4}--\eqref{MOc9}, as the following linear inequalities:
\begin{align}
& \exists \ \varepsilon^{\tinyL}, \varepsilon^{\tinyU}, S, s_0, R: \nonumber \\
& Re \leq H\varepsilon^{\tinyL} - Tx - Ws_0, \label{ref-note-5} \\
& R \geq WS - HE, \ R \geq 0. \label{ref-note-6}
\end{align}
\end{subequations}

Second, we recast the adjustable DR joint chance constraint \eqref{Pc2} as second-order conic constraints. We present this result in the following theorem and its proof in Appendix \ref{apx:theorem-1}.
\begin{theorem} \label{theorem-1}
If $u > 2/3$, then, for all $t \in [T]$, chance constraint \eqref{Pc2} is equivalent to the following constraints:\\
\begin{subequations}
\begin{align}
& \left\| \left[ \begin{matrix} \sqrt{\tfrac{8}{3}}\\ r_{kt}-z_{kt}
\end{matrix} \right] \right\|_2 \leq r_{kt}+z_{kt}, \ \forall k \in [K], \label{The11}\\
& \left\| \left[ \begin{matrix} s_{kt}-1\\ 2z_{kt}
\end{matrix} \right] \right\|_2 \leq s_{kt}+1, \ \forall k \in [K], \label{The12}\\
& \sigma_{kt} r_{kt} \leq \mu_{kt}-\varepsilon_{kt}^{\tinyL}, \ \forall k \in [K], \label{The13}\\
& \sigma_{kt} r_{kt} \leq \varepsilon_{kt}^{\tinyU}-\mu_{kt}, \ \forall k \in [K], \label{The14}\\
& \sum_{k \in [K]} s_{kt} \leq 1-u, \label{The15}\\
& r_{kt},s_{kt},z_{kt} \geq 0, \ \forall k \in [K]. \label{The16}
\end{align}
\end{subequations}
\end{theorem}

To summarize, the DRCO model \eqref{drco-obj}--\eqref{drco-con} is equivalent to the following second-order conic program:
\begin{subequations}
\begin{align}
\min_{\substack{\hat{p},S,s_0,R,u,\\ \varepsilon^{\tinyL},\varepsilon^{\tinyU},r,s,z}} & \ \sum_{t \in [T]} \sum_{i \in [I]} C_i(\hat{p}_{it}) - \delta u \label{drco-ref-obj} \\
\text{s.t.} \ & \ \mbox{\eqref{EDc1}--\eqref{EDc4}}, \ \mbox{\eqref{ref-note-5}--\eqref{ref-note-6}}, \ \eqref{MOc10}, \ \mbox{\eqref{Pc3}}, \ \mbox{\eqref{The11}--\eqref{The16}}. \label{drco-ref-con}
\end{align}
\end{subequations}

\section{Extension to Alternative Operational Risks} \label{sec:dr-expected-cost}
We extend the DRCO model \eqref{drco-obj}--\eqref{drco-con} by considering alternative operational risks of the chance constraint \eqref{Pc2}, which computes the expected costs incurred by overestimating/underestimating the renewable energy. We note that such operational risks are studied in~\cite{wang2016risk}. In this paper, we study the DR counterpart of the risks based on the ambiguity set $\mathcal{D}$ defined in \eqref{ambiguity-set}. Specifically, the DR expected cost of overestimation/underestimation are defined as
\begin{align*}
P^+(\varepsilon^{\tinyL}, \varepsilon^{\tinyU}) := \ \sup_{\mathbb{P} \in \mathcal{D}} \mathbb{E_P} \left[ \sum_{k \in [K]} \sum_{t \in [T]} c^+_{kt} \left[ \varepsilon_{kt}^{\tinyL} - \varepsilon_{kt} \right]^+ \right], \\ P^-(\varepsilon^{\tinyL}, \varepsilon^{\tinyU}) := \ \sup_{\mathbb{P} \in \mathcal{D}} \mathbb{E_P} \left[ \sum_{k \in [K]} \sum_{t \in [T]} c^-_{kt} \left[ \varepsilon_{kt} - \varepsilon_{kt}^{\tinyU} \right]^+ \right],
\end{align*}
where $[x]^+ = \max\{x, 0\}$ for $x \in \mathbb{R}$. When $\varepsilon_{kt} \notin [\varepsilon_{kt}^{\tinyL}, \varepsilon_{kt}^{\tinyU}]$, emergency regulations (e.g., renewable generation curtailment, fast-starting units, and load shedding) may be needed to recover the operational feasibility. Accordingly, the cost coefficients $c^+_{kt}$ and $c^-_{kt}$ should be estimated based on the corresponding regulations (e.g., the opportunity/penalty cost of curtailing renewable generation, the estimated real-time price of using fast-starting units, and the penalty cost of shedding load). Then, the DRCO model \eqref{drco-obj}--\eqref{drco-con} can be extended by incorporating $P^{\pm}(\varepsilon^{\tinyL}, \varepsilon^{\tinyU})$ as follows: %As in the chance constraint \eqref{Pc2}, we employ $\mathcal{D}$ to model the ambiguity of $\mathbb{P}$ in the definition of $P^+(\varepsilon^{\tinyL}, \varepsilon^{\tinyU})$ and $P^-(\varepsilon^{\tinyL}, \varepsilon^{\tinyU})$.
\begin{align*}
\min_{\substack{\hat{p},B,b,\\ \varepsilon^{\tinyL},\varepsilon^{\tinyU},u}} & \ \sum_{t \in [T]} \sum_{i \in [I]} C_i(\hat{p}_{it}) - \delta u + \delta^+ P^+(\varepsilon^{\tinyL}, \varepsilon^{\tinyU}) + \delta^- P^-(\varepsilon^{\tinyL}, \varepsilon^{\tinyU}) \\
\text{s.t.} \ & \ \mbox{\eqref{EDc1}--\eqref{EDc4}}, \ \mbox{\eqref{MOc4}--\eqref{MOc9}}, \ \eqref{MOc10}, \ \mbox{\eqref{Pc2}--\eqref{Pc3}},
\end{align*}
where $\delta^+$ and $\delta^-$ represent the weights on the expected costs of overestimation/underestimation, respectively.

We compute $P^+(\varepsilon^{\tinyL}, \varepsilon^{\tinyU})$ and $P^-(\varepsilon^{\tinyL}, \varepsilon^{\tinyU})$ by solving conic programs. We present this result in the following theorem and its proof in Appendix \ref{apx:theorem-2}. Accordingly, the extended DRCO model presented above can be recast as a conic program that can be efficiently solved by commercial solvers.
%We compute $P^+(\varepsilon^{\tinyL}, \varepsilon^{\tinyU})$ and $P^-(\varepsilon^{\tinyL}, \varepsilon^{\tinyU})$ by solving conic programs. We present this result in the following theorem and omit its proof due to the space limit\footnote{Please see https://arxiv.org/ftp/arxiv/papers/xxxx/xxxx.xxxxx.pdf for a detailed proof.}. Accordingly, the extended DRCO model presented above can be recast as a conic program that can be efficiently solved by commercial solvers.
\begin{theorem} \label{theorem-2}
Let $g(\{(c_{kt}, \tau_{kt})\}_{k \in [K], t \in [T]})$ represent the optimal value of the following conic program:
\begin{subequations}
\begin{align}
\min_{\substack{\pi_{\ell kt}, \Lambda_{ktij}}} & \ \sqrt{3}\sum_{k \in [K]} \sum_{t \in [T]} c_{kt} \sigma_{kt} (\pi_{1kt} + \pi_{3kt} + 1) \\
\mbox{s.t.} & \ \Lambda_{kt00} = \tau_{kt}, \ \forall k \in [K], \ \forall t \in [T], \\
& \ \left\| \begin{bmatrix} \pi_{2kt} \\ \pi_{1kt} - \pi_{3kt} + 1 \end{bmatrix} \right\|_2 \leq \pi_{1kt} + \pi_{3kt} + 1, \nonumber \\
& \ \forall k \in [K], \ \forall t \in [T], \\
& \ \sum_{i,j: \ i+j = 2\ell-1} \Lambda_{ktij} = 0, \nonumber \\
& \ \forall \ell = 1, 2, 3, \ \forall k \in [K], \ \forall t \in [T], \\
& \ \sum_{i,j: \ i+j = 2\ell} \Lambda_{ktij} = \pi_{\ell kt}, \nonumber \\
& \ \forall \ell = 1, 2, 3, \ \forall k \in [K], \ \forall t \in [T], \\
& \ \Lambda_{kt} \in \mathbb{S}^{4 \times 4}_+, \ \forall k \in [K], \ \forall t \in [T],
\end{align}
\end{subequations}
where $\mathbb{S}^{4 \times 4}_+$ represents the cone of all $4\times 4$ positive semidefinite matrices. Then, we have
\begin{align*}
& P^+(\varepsilon^{\tinyL}, \varepsilon^{\tinyU}) = g\left(\left\{\left(c^+_{kt}, \frac{\mu_{kt} - \varepsilon^{\tinyL}_{kt}}{\sqrt{3} \sigma_{kt}}\right)\right\}_{k \in [K], t \in [T]}\right), \\
& P^-(\varepsilon^{\tinyL}, \varepsilon^{\tinyU}) = g\left(\left\{\left(c^-_{kt}, \frac{\varepsilon^{\tinyU}_{kt} - \mu_{kt}}{\sqrt{3} \sigma_{kt}}\right)\right\}_{k \in [K], t \in [T]}\right).
\end{align*}
Furthermore, $P^{\pm}(\varepsilon^{\tinyL}, \varepsilon^{\tinyU})$ can be conservatively approximated by piecewise linear functions of $(\varepsilon^{\tinyL}, \varepsilon^{\tinyU})$ with arbitrary precision.
\end{theorem}

%%%%%%%%%%%%%%%%%%%%%%%%%%%%%%%%%%%%%%%%%%%%%%%%%%%%%%
%%%%%%%%%%%%%%%%%%%%%%%%%%%%%%%%%%%%%%%%%%%%%%%%%%%%%%
\section{Case Studies} \label{sec:computation}
%We carry out numerical case studies on modified IEEE 14-bus and IEEE 118-bus systems to demonstrate the effectiveness of the proposed approach. All programs are developed using MATLAB2014a and solved by MOSEK via YAlMIP 11.5 on a laptop with a 2.7GHz Intel Core i5 CPU and 8GB RAM.
We carry out numerical case studies on modified IEEE 14-bus and IEEE 118-bus systems. All programs are developed using MATLAB2014a and solved by MOSEK via YAlMIP 11.5 on a laptop with a 2.7GHz Intel Core i5 CPU and 8GB RAM.

\subsection{The Modified IEEE 14-bus System}
In this system, there are 20 transmission lines and 5 generators (G1--G5) providing corrective power re-dispatch. The generators and network characteristics can be found in MATPOWER~\cite{zimmerman2011matpower}. Two wind power farms with 80 MW (W1) and 100MW (W2) installed capacity are connected to the system at nodes 5 and 7, respectively. The load profile is from~\cite{wang2017robustrcuc} and scaled by a factor of 0.1. We set $\Delta_{d} = 60\mbox{min}$, $\Delta_{r} = 5\mbox{min}$, $T = 24\mbox{hr}$, and $p^{\text{min}}_i = 0.1p^{\text{max}}_i$ for all $i \in [I]$. %We assume that at initial time 0, operators look 24 dispatch periods ahead and decide the dispatch strategies including the outputs of generators and the admissible ranges of wind power.
%%%%%%%%%%%%%%%%%%%%%

\subsection{The Cost-Utilization Frontier and the DNE Limits}
To demonstrate the trade-off between the dispatch cost and the utilization of renewable energy, we generate a cost-utilization frontier by gradually increasing the value of $\delta$ and re-solving the DRCO model \eqref{drco-obj}--\eqref{drco-con} for each $\delta$. To this end, we first obtain the wind power forecast of W1 and W2 from the NREL Eastern Wind Dataset~\cite{nreleasternwind}. We generate a set of wind power prediction error data by using Gaussian distribution, whose mean is set to be $0$ for all $t \in [T]$ and variance increases from 10\% of the installed capacity by 0.1\% as $t$ increases from $1$ to $T$. Second, we divide the data into two parts. We use the first part to calibrate the ambiguity set $\mathcal{D}$ based on the empirical mean and variance. Then, for fixed $\delta$, we solve the DRCO model to obtain the optimal DNE limits $[\varepsilon^{\tinyL*}, \varepsilon^{\tinyU*}]$ and the minimum power dispatch cost. We use the second part of the data to obtain an out-of-sample empirical estimate of the renewable utilization probability $\mathbb{P}\{\varepsilon \in [\varepsilon^{\tinyL*}, \varepsilon^{\tinyU*}]\}$. Third, we repeat the second step by gradually increasing the value of $\delta$ from $1$ to $38000$. We set the step length as 100 when $\delta \leq 1000$, $400$ when $\delta \in (1000, 5000]$, $1000$ when $\delta \in (5000, 10000]$, and $4000$ when $\delta > 10000$. Accordingly, we obtain 33 groups of minimum power dispatch costs, optimal DNE limits, and the corresponding renewable utilization probabilities.

\begin{figure}[h]
\centering
\includegraphics[width=0.8\linewidth]{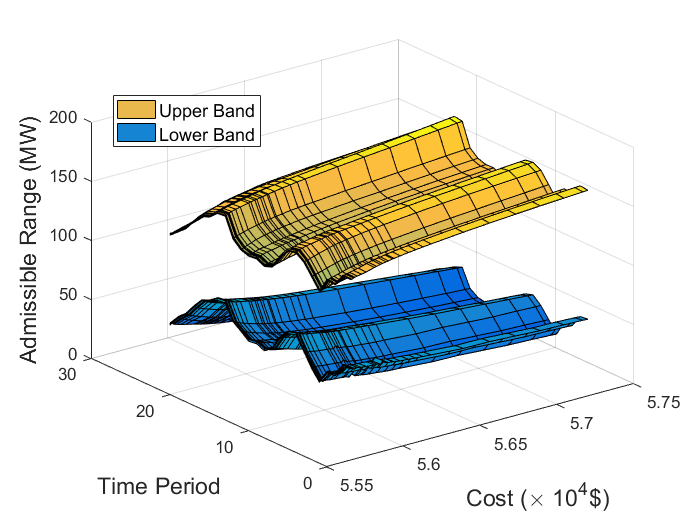}
\caption{Minimum Dispatch Costs vs. Admissible Ranges of Wind Power, with $r^{\text{up}}_i = r^{\text{dn}}_i = 2.0\% \times p^{\text{max}}_i$ for all $i \in [I]$}
\label{fig14range}
\end{figure}

In Fig.~\ref{fig14range}, we display the minimum power dispatch cost and the optimal DNE limits under various $\delta$ values. For intuitive presentation, we shift the DNE limits to obtain the admissible ranges of wind power $[w^{\tinyL}_t, w^{\tinyU}_t] := [\sum_{k \in [K]} (\mu_{kt} + \varepsilon^{\tinyL}_{kt}), \sum_{k \in [K]} (\mu_{kt} + \varepsilon^{\tinyU}_{kt})]$. The interpretation of each pair of points $(t, C_{\text{min}}, w^{\tinyL}_t)$ and $(t, C_{\text{min}}, w^{\tinyU}_t)$ is that, during time period $t$, we need to spend at least $C_{\text{min}}$ on power dispatch in order to accommodate any total wind power output within the interval $[w^{\tinyL}_t, w^{\tinyU}_t]$. From Fig.~\ref{fig14range}, we observe that the admissible range of wind power broadens as the minimum dispatch cost increases. This indicates that the power system can become more flexible as we invest more on power dispatch.

%%%%%%%%%%%%%%%%%%%%%%%%

\begin{figure}[h]
\centering
\includegraphics[width=0.8\linewidth]{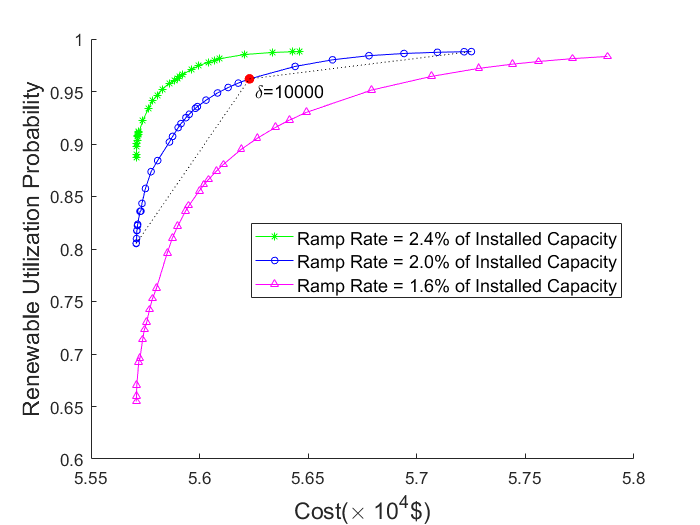}
\caption{The Cost-Utilization Frontier under Various Ramping Capabilities}
\label{fig14con}
\end{figure}
%%%%%%%%%%%%%%%%%%%%%%%%

To take a closer look on the trade-off between the dispatch cost and the renewable utilization, we display the cost-utilization frontiers under various ramping capabilities in Fig.~\ref{fig14con}. Specifically, we consider three ramping capabilities in which $r^{\text{up}}_i$ and $r^{\text{dn}}_i$ equal 1.6\%, 2.0\%, and 2.4\% of $p^{\text{max}}_i$ for all $i \in [I]$, respectively. For each capability and for each $\delta$ value, we depict the minimum power dispatch cost versus the smallest renewable utilization probability in all time periods. From Fig.~\ref{fig14con}, we first observe that the renewable utilization increases as the dispatch cost increases, confirming our observation from Fig.~\ref{fig14range}. Second, the increasing trend of renewable utilization diminishes as the dispatch cost increases. Take the middle curve with 2\% ramping capability for example. On this curve, we highlight two segments with $\delta \leq 10,000$ and $\delta > 10,000$, respectively. The first segment reflects a 19.5\% increase in renewable utilization with only an 0.9\% increase in the power dispatch cost (i.e., increasing by \$521). This translates into a 0.037\%/\$ increasing rate of the renewable utilization. On the contrary, the second segment reflects a 2.7\% increase in renewable utilization with a 1.8\% increase in the power dispatch cost (i.e., increasing by \$1021). This translates into a 0.003\%/\$ increasing rate of the renewable utilization. This observation indicates that a small additional investment on power dispatch can quickly enhance the renewable utilization, but this investment becomes less efficient when the utilization is already high. Third, we observe from Fig.~\ref{fig14con} that the frontier rises as the ramping capability increases. For example, to achieve a 95\% renewable utilization, it costs $5.67\times 10^4\$$ when the ramping capability is 1.6\%, $5.61\times 10^4\$$ when the ramping capability is 2.0\%, and $5.58\times 10^4\$$ when the ramping capability is 2.4\%. This observation indicates that a small enhancement on the ramping capability can significantly improve the cost-effectiveness of utilizing renewable energy.

\begin{figure}[h]
        %\centering
        \begin{subfigure}[b]{0.255\textwidth}
                \includegraphics[width=\textwidth]{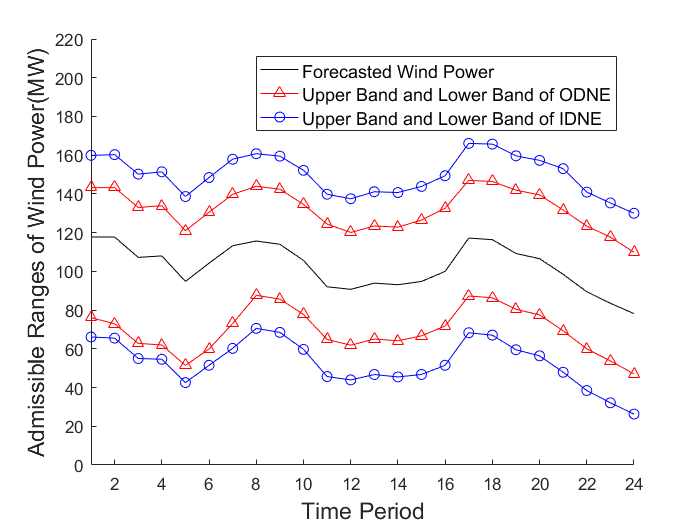}
        \end{subfigure} \hspace{-0.8cm} %
        ~ %add desired spacing between images, e. g. ~, \quad, \qquad etc.
        \begin{subfigure}[b]{0.255\textwidth}
                \includegraphics[width=\textwidth]{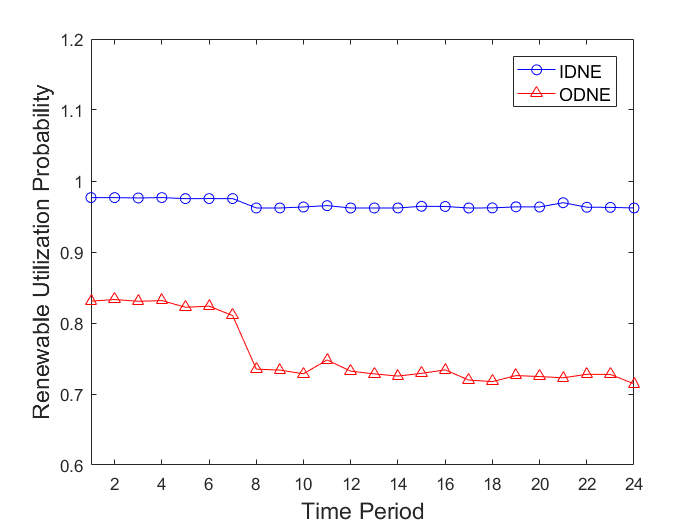}
        \end{subfigure}
        \caption{Comparison on (a) Admissible Range and (b) Renewable Utilization Probability.}
        \label{fig14costandcon}
\end{figure}

%%%%%%%%%%%%%%%%%%%%%%%%%%%%%%%%%
%%%%%%%%%%%%%%%%%%%%%%%%
\subsection{Comparisons with the Original DNE Limit Approach}
We compare the proposed DRCO model (termed the IDNE approach) with the original DNE approach (termed ODNE) in~\cite{zhao2015dnelimit}, which computes the DNE limits based on a given dispatch strategy without explicitly modeling the renewable ambiguity. We randomly generate 5,000 out-of-sample scenarios of wind prediction errors from the hypothetical Gaussian distribution and compare (i) the optimal DNE limits, (ii) the renewable utilization probability, and (iii) the actual cost incurred in each scenario. The actual cost consists of the pre-dispatch cost, the corrective re-dispatch cost, and the penalty costs which are incurred (a) when $\varepsilon_{kt} < \varepsilon^{\tinyL}_{kt}$, the load shedding takes place at a cost of 2,000\$/MW and (b) when $\varepsilon_{kt} > \varepsilon^{\tinyU}_{kt}$, the renewable energy is curtailed at a cost of 100\$/MW. In all comparisons, we set $\delta = 10,000$ and $r^{\text{up}}_i=r^{\text{dn}}_i=2.0\% \times p^{\text{max}}_i$ for all $i \in [I]$.

First, we compare the admissible ranges of wind power and the corresponding out-of-sample renewable utilization probabilities in Fig.~\ref{fig14costandcon}. From this figure, we observe that the proposed IDNE approach yields wider admissible ranges, and so higher renewable utilization, than ODNE does. For example, IDNE can consistently accommodates more than 95\% of the wind power throughout the 24 time periods, while ODNE accommodates less than 90\% and shows a decreasing trend in renewable utilization as $t$ increases.
%We adopt $\delta$ as 10000 for the proposed integrated do-not-exceed (IDNE) limit approach as the red point in Fig.\ref{fig14con} when the confidence level of 96.33\% can be ensured. The original do-not-limit (ODNE) approach for comparison refers to \cite{zhao2015dnelimit}, which gives optimal accessible wind power variation ranges based on economic dispatch strategies without considering uncertainty.
%%%%%%%%%%%%%%%%%%%%%%%%%%%%%%%%%
\begin{figure}[h]
\centering
\includegraphics[width=0.8\linewidth]{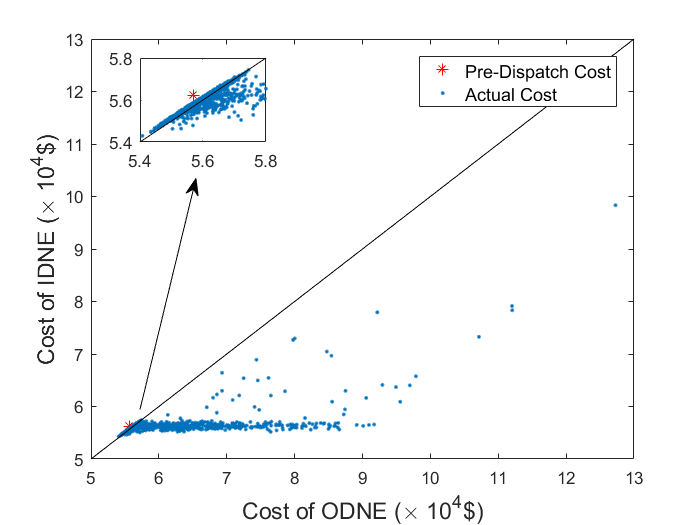}
\caption{Comparison on Actual Cost}
\label{fig14cost}
\end{figure}

\begin{table}[h]
	\centering
	\caption{Comparisons of two approaches}
	\begin{tabular}{ccccc}
	\toprule
	 Approach & AvgC (\$) & MaxC (\$) & AvgLS (MW) & AvgWC (MW) \\
	\midrule
     IDNE & 56,004 & 98,333 &  0.0349 &  0.0276  \\
     ODNE & 57,030 & 127,277 &  0.5686 & 0.2389  \\
     %\hline
     %Difference & -672 & -25092 & -0.3716 &-0.1071\\
	\bottomrule
	\end{tabular}
	\label{Table1}
\end{table}
%%%%%%%%%%%%%%%%%%%%%%%%%%%%%%%
Second, we compare the actual cost in Fig.~\ref{fig14cost} and Table~\ref{Table1}. In Fig.~\ref{fig14cost}, we plot the IDNE actual cost versus the ODNE actual cost for all 5,000 scenarios and the 45-degree reference line represents that the two costs agree. From this figure, we observe that most points distribute around or below the reference line, indicating that IDNE is likely to outperform ODNE in out-of-sample tests. In addition, most points line up along the horizontal line of $5.6\times 10^4\$$, i.e., the IDNE pre-dispatch cost. This indicates that the IDNE yields stable actual costs with small variations. That is, the proposed DR approach provide stable and predictable out-of-sample performance. In Table~\ref{Table1}, we report the average actual cost (AvgC), the maximum actual cost (MaxC), the average load shedding (AvgLS), and the average wind curtailment (AvgWC) among the 5,000 scenarios. This table confirms the observations on the actual cost from Fig.~\ref{fig14cost}, and further demonstrates that IDNE incurs one order of magnitude less load shedding as well as wind curtailment than ODNE does.

%%%%%%%%%%%%%%%%%%%%%%%%%%%%%%%
\begin{figure}[h]
\centering
\includegraphics[width=0.7\linewidth]{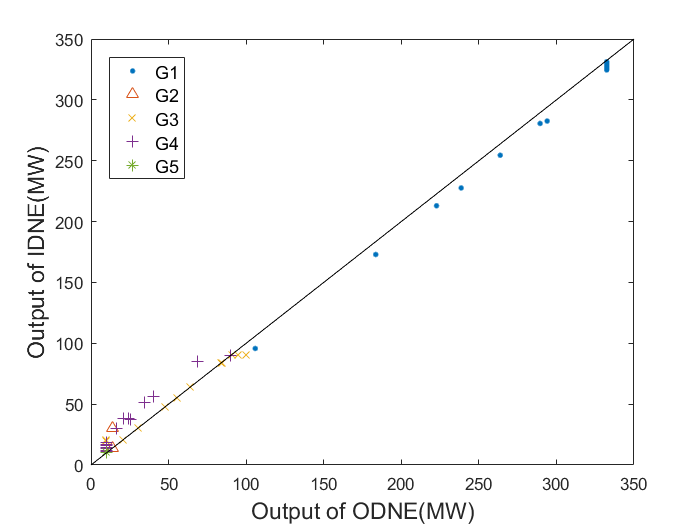}
\caption{Comparison on the Generation Amounts of the Thermal Units}
\label{fig14strategies}
\end{figure}
%%%%%%%%%%%%%%%%%%%%%%%%%%%%%%%
Third, we compare the optimal pre-dispatch strategies of IDNE and ODNE in Fig.~\ref{fig14strategies}, from which we observe that the pre-dispatch generation amounts of G1 and G3 under IDNE are lower than those under ODNE in most time periods. This strategy enhances the system flexibility under IDNE by preserving more ramping capability, especially when the load is high. Take the time period $t = 11$ for example, in which the load is high and G1, G3 under ODNE reach their maximum generation capacities. In this case, the upward ramping capability of ODNE becomes scarce. On the contrary, IDNE sets a lower pre-dispatch generation amounts of G1 and G3, and so preserves more (upward) ramping capability. This demonstrates how the power dispatch and DNE limits can coordinate in the proposed DRCO model to enhance the system flexibility. Finally, it takes 30.80 CPU seconds on average and 33.42 CPU seconds at maximum to solve the DRCO model with various values of $\delta$, which verifies the tractability of the proposed approach.

%Applying the proposed IDNE limit approach and the ODNE approach respectively, the dispatch strategies of 5 generators are displayed in Fig.\ref{fig14strategies}. During dispatch periods except for period 14 and 17, the outputs of G1 under IDNE limit approach are lower than those under ODNE limit approach, and simultaneously the outputs of G3 or G4 are higher under IDNE limit approach than under the ODNE limit approach. Such difference in dispatch strategies increases the flexibility of the system through preserving more adjustable capacities. As a result, the proposed IDNE limit approach broadens the admissible ranges of wind power variation and improves the confidence level when compared to the ODNE limit approach, as shown in Fig.\ref{fig14costandcon}. But as a sacrifice, the dispatch cost of the IDNE limit approach is \$56355 which is bigger than that of ODNE limit approach \$55854.

%%%%%%%%%%%%%%%%%%%%%%%%%%%%%%%%%
\begin{figure}[h]
\centering
\includegraphics[width=0.7\linewidth]{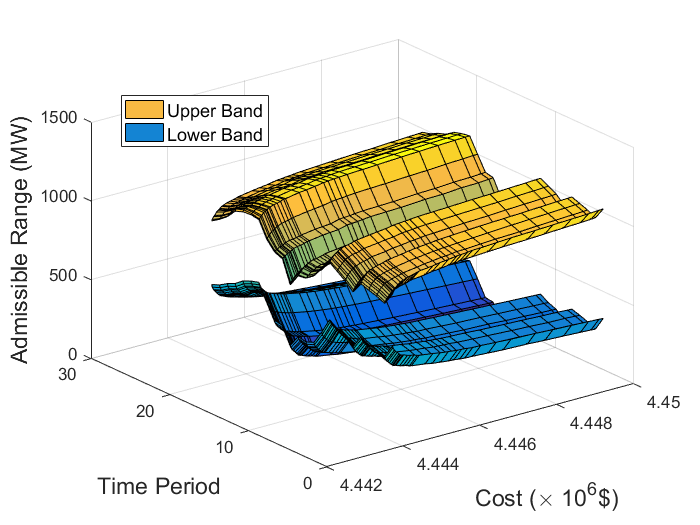}
\caption{Minimum Dispatch Cost vs. Admissible Ranges of Wind Power based on the Modified IEEE-118 System}
\label{fig118range}
\end{figure}

\begin{figure}[h]
\centering
\includegraphics[width=0.7\linewidth]{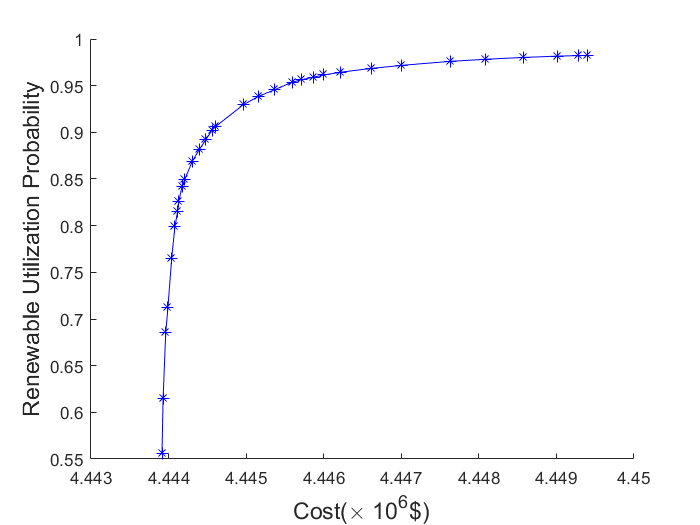}
\caption{The Cost-Utilization Frontier based on the Modified IEEE-118 System}
\label{fig118cost}
\end{figure}
%%%%%%%%%%%%%%%%%%%%%%%%%%%%%%%%%

%To further illustrate the advantages of the proposed IDNE limit method, we generate 5,000 realized wind power output curves based on the hypothetical Gaussian distributions and calculate their actual costs taking the load shedding and wind power curtailment penalty costs into account. The penalty price for load shedding and wind power curtailment are set 2000\$/MWh and 100\$/MWh respectively. For each realized wind power curve, the corrective dispatches under two approaches are performed and the actual costs are evaluated. Fig.\ref{fig14cost} displays the actual costs of two approaches under 5,000 simulations. It can be seen that the actual costs under the strategy of the proposed IDNE limit approach presents more robust and smoothing than that of the ODNE method which has many spikes. The average actual cost (AvgC), the maximum actual cost (MaxC), the average amount of load shedding (AvgLS) and the average amount of wind curtailment(AvgWC) of 5,000 times tests under two approaches are calculated as shown in Tab.\ref{Table1}. the proposed method is \$673 smaller than that of the ODNE limit method as And the average load shedding and wind curtailment amount of under the proposed method are also smaller than those of ODNE limit method.

\subsection{The Modified IEEE 118-bus System}
In the modified IEEE 118-bus system, there are 186 transmission lines and 54 generators providing corrective re-dispatch. Three wind farms with an identical 300 MW installed capacity are connected to the system at nodes 18, 32, and 88, respectively. The generator and network characteristics are from MATPOWER 5.1~\cite{zimmerman2011matpower} and the load profile is the same as in~\cite{wang2017robustrcuc}. In addition, the wind forecasts are from the NREL Eastern Wind Dataset~\cite{nreleasternwind} and the mean and variance of the wind power prediction error are set as in the previous case study. We display the minimum power dispatch cost and the optimal DNE limits under various $\delta$ values in Fig.~\ref{fig118range} and the cost-utilization frontier in Fig.~\ref{fig118cost}. We make similar observations from these two figures on the trade-off between the power dispatch cost and the renewable utilization. Finally, it takes 296.38 CPU seconds on average and 307.51 CPU seconds at maximum to solve the DRCO model with various values of $\delta$, which verifies the scalability of the proposed approach.

% An example of a double column floating figure using two subfigures.
% (The subfig.sty package must be loaded for this to work.)
% The subfigure \label commands are set within each subfloat command,
% and the \label for the overall figure must come after \caption.
% \hfil is used as a separator to get equal spacing.
% Watch out that the combined width of all the subfigures on a
% line do not exceed the text width or a line break will occur.
%
%\begin{figure*}[!t]
%\centering
%\subfloat[Case I]{\includegraphics[width=2.5in]{box}%
%\hfil
%\subfloat[Case II]{\includegraphics[width=2.5in]{box}%
%\label{fig_second_case}}
%\caption{Simulation results for the network.}
%\label{fig_sim}
%\end{figure*}
%
% Note that often IEEE papers with subfigures do not employ subfigure
% captions (using the optional argument to \subfloat[]), but instead will
% reference/describe all of them (a), (b), etc., within the main caption.
% Be aware that for subfig.sty to generate the (a), (b), etc., subfigure
% labels, the optional argument to \subfloat must be present. If a
% subcaption is not desired, just leave its contents blank,
% e.g., \subfloat[].

%%%%%%%%%%%%%
\section{Conclusion and Future Research} \label{sec:conclusion}
We propose a DRCO model for power dispatch and DNE limits. Our model incorporates an adjustable DR joint chance constraint to explicitly measure the utilization of renewable energy. By using ADR, we derive a second order conic program that conservatively approximates the DRCO model. The case studies based on modified IEEE 14-bus and 118-bus systems demonstrate the effectiveness and computational tractability of the proposed approach. Future research includes alternative ambiguity sets and the corresponding DRCO models.

%%%%%%%%%%%%%
\appendices
\section{Proof of Theorem \ref{theorem-1}} \label{apx:theorem-1}
First, we observe that ambiguity set $\mathcal{D}$ satisfies Assumption (A1) in~\cite{xie2017optimized}. Hence, by Theorem 3 in~\cite{xie2017optimized}, the chance constraint \eqref{Pc2} is equivalent to its Bonferroni approximation:
\begin{subequations}
\begin{align}
& \inf_{\mathbb{P}_{kt} \in \mathcal{D}_{kt}} \mathbb{P}_{kt} (\varepsilon_{kt} \in [\varepsilon_{kt}^{\tinyL},\varepsilon_{kt}^{\tinyU}] ) \geq 1-s_{kt}, \nonumber \\
& \forall k \in [K], \ \forall t \in [T], \label{Bapp1}\\
& \sum_{k \in [K]} s_{kt} \leq 1-u, \ \forall t \in [T], \label{Bapp2}\\
& s_{kt} \geq 0, \ \forall k \in [K], \ \forall t \in [T], \label{Bapp3}
\end{align}
where $\mathbb{P}_{kt}$ represents the (marginal) probability distribution of each $\varepsilon_{kt}$ and $\mathcal{D}_{kt} = \{\mathbb{P}_{kt}: \mathbb{E}_{\mathbb{P}_{kt}}[\varepsilon_{kt}] = \mu_{kt}, \mbox{Var}(\varepsilon_{kt}) = \sigma_{kt}^2,$ $\mathbb{P}_{kt} \mbox{ is unimodal about $\mu_{kt}$} \}$.

Second, chance constraint \eqref{Bapp1} is equivalent to:
\begin{align}
\inf_{\mathbb{P}_{kt} \in \mathcal{D}_{kt}} \mathbb{P}_{kt} & \left( |\varepsilon_{kt}-\mu_{kt}|  \leq \min \{\mu_{kt}-\varepsilon_{kt}^{\tinyL},\varepsilon_{kt}^{\tinyU}-\mu_{kt}\}  \right) \notag \\
&\geq 1-s_{kt}, \ \forall k \in [K], \ \forall t \in [T], \label{Bapp1app}
\end{align}
where $1-s_{kt} \geq u > 2/3$ due to \eqref{Bapp2}. By the Gauss inequality~\cite{gauss1823theoria}, we recast \eqref{Bapp1app} as
\begin{align}
1-\frac{4}{9\lambda_{kt}^2} \geq 1-s_{kt}, \ \forall k \in [K], \ \forall t \in [T], \label{ka1}
\end{align}
where $\lambda_{kt} := \min \{\mu_{kt}-\varepsilon_{kt}^{\tinyL},\varepsilon_{kt}^{\tinyU}-\mu_{kt}\}/\sigma_{kt}$.
\end{subequations}

Third, we recast inequality \eqref{ka1} as second-order conic constraints \eqref{The11}--\eqref{The14} by introducing auxiliary variables $r_{kt}$ and $z_{kt}$ (see~\cite{ben2001lectures}).

\section{Proof of Theorem \ref{theorem-2}} \label{apx:theorem-2}
First, for given $\tau \in \mathbb{R}_+$, we compute the worst-case expectation $J(\tau) := \sup_{\mathbb{P} \in \mathcal{D}_0}\mathbb{E_P}[\varepsilon - \tau]^+$ with $\mathcal{D}_0 = \{\mathbb{P}: \mathbb{E_P}[\varepsilon] = 0, \mbox{Var}(\varepsilon) = 1/3, \varepsilon \mbox{ is unimodal about } 0 \}$. To this end, by the unimodality of $\varepsilon$, there exists a random variable $\zeta \in \mathbb{R}$ such that $\varepsilon = U\zeta$, where $U$ is uniform on $(0, 1)$ and independent of $\zeta$ (see~\cite{dharmadhikari1988unimodality}). It follows that $\mathbb{E_P}[\zeta] = 0$, $\mbox{Var}(\zeta) = 1$, and $\mathbb{E_P}[\varepsilon - \tau]^+ = \mathbb{E_P}[h(\zeta)]$, where $$
h(\zeta) = \left\{\begin{array}{ll} 0 & \mbox{if $\zeta \leq 0$,}\\
\left[1 - \frac{\tau}{\zeta}\right]^+ & \mbox{if $\zeta > 0$.} \end{array}\right.
$$
We compute $\sup_{\mathbb{P} \in \mathcal{D}_0}\mathbb{E_P}[h(\zeta)]$ by formulating the following optimization problem:
\begin{subequations}
\begin{align}
\max_{\mathbb{P}} & \ \int_{\mathbb{R}} h(\zeta) d\mathbb{P} \label{apx-thm-2-note-1} \\
\mbox{s.t.} & \ \int_{\mathbb{R}} \zeta^2 d\mathbb{P} = 1, \label{apx-thm-2-note-2} \\
& \ \int_{\mathbb{R}} \zeta d\mathbb{P} = 0, \label{apx-thm-2-note-3} \\
& \ \int_{\mathbb{R}} d\mathbb{P} = 1, \label{apx-thm-2-note-4}
\end{align}
whose dual formulation is
\begin{align}
\min_{\pi} & \ \pi_3 + \pi_1 + 1 \label{apx-thm-2-note-5} \\
\mbox{s.t.} & \ \pi_3 \zeta^2 + \pi_2 \zeta + \pi_1 + 1 \geq 0, \ \forall \zeta \geq 0, \label{apx-thm-2-note-6} \\
& \ \pi_3 \zeta^2 + \pi_2 \zeta + \pi_1 + 1 \geq 1 - \frac{\tau}{\zeta}, \ \forall \zeta \geq 0, \label{apx-thm-2-note-7}
\end{align}
where dual variables $\pi_3$, $\pi_2$, and $\pi_1 + 1$ are associated with primal constraints \eqref{apx-thm-2-note-2}--\eqref{apx-thm-2-note-4}, respectively. Dual constraint \eqref{apx-thm-2-note-6} is equivalent to $\pi_1 + 1 - (\pi_2)^2/(4\pi_3) \geq 0$, which is further equivalent to
\begin{equation}
\left\| \begin{bmatrix} \pi_2 \\ \pi_1 - \pi_3 + 1 \end{bmatrix} \right\|_2 \leq \pi_1 + \pi_3 + 1. \label{apx-thm-2-note-8}
\end{equation}
In addition, dual constraint \eqref{apx-thm-2-note-7} is equivalent to $\pi_3 \zeta^3 + \pi_2 \zeta^2 + \pi_1 \zeta + \tau \geq 0$ for all $\zeta \geq 0$, which is further equivalent to
\begin{align}
& \ \exists \Lambda \in \mathbb{S}^{4 \times 4}_+ \mbox{ such that:} \nonumber \\
& \ \Lambda_{00} = \tau, \label{apx-thm-2-note-9} \\
& \ \sum_{i,j: \ i+j = 2\ell-1} \Lambda_{ij} = 0, \ \forall \ell = 1, 2, 3, \label{apx-thm-2-note-10} \\
& \ \sum_{i,j: \ i+j = 2\ell} \Lambda_{ij} = \pi_{\ell}, \ \forall \ell = 1, 2, 3 \label{apx-thm-2-note-11}
\end{align}
by Proposition 3.1(b) in~\cite{bertsimas2005optimal}. It follows that $J(\tau) \equiv \sup_{\mathbb{P} \in \mathcal{D}_0}\mathbb{E_P}[h(\zeta)]$ equals the optimal value of the following conic program:
\begin{align}
\min_{\pi, \Lambda} \ & \ \pi_3 + \pi_1 + 1 \label{apx-thm-2-note-12} \\
\mbox{s.t.} \ & \ \mbox{\eqref{apx-thm-2-note-8}--\eqref{apx-thm-2-note-11}}, \ \Lambda \in \mathbb{S}^{4 \times 4}_+. \label{apx-thm-2-note-13}
\end{align}

Second, we have $P^-(\varepsilon^{\tinyL}, \varepsilon^{\tinyU}) = \sum_{k \in [K]}\sum_{t \in [T]} c^-_{kt} \sup_{\mathbb{P} \in \mathcal{D}}$ $\mathbb{E_P}[\epsilon_{kt} - \epsilon^{\tinyU}_{kt}]^+$ because $\mathcal{D}$ is separable over indices $k$ and $t$. But
\begin{align*}
& \sup_{\mathbb{P} \in \mathcal{D}}\mathbb{E_P}\left[\epsilon_{kt} - \epsilon^{\tinyU}_{kt}\right]^+ \\
= \ & \sqrt{3}\sigma_{kt} \sup_{\mathbb{P} \in \mathcal{D}}\mathbb{E_P}\left[\left(\frac{\epsilon_{kt} - \mu_{kt}}{\sqrt{3}\sigma_{kt}}\right) - \left(\frac{\epsilon^{\tinyU}_{kt} - \mu_{kt}}{\sqrt{3}\sigma_{kt}}\right)\right]^+ \\
= \ & \sqrt{3}\sigma_{kt} \ J\left(\frac{\epsilon^{\tinyU}_{kt} - \mu_{kt}}{\sqrt{3}\sigma_{kt}}\right)
\end{align*}
because random variable $(\epsilon_{kt} - \mu_{kt})/(\sqrt{3}\sigma_{kt})$ has mean $0$, variance $1/3$, and is unimodal about $0$. Hence,
\begin{align*}
P^-(\varepsilon^{\tinyL}, \varepsilon^{\tinyU}) = & \sqrt{3} \sum_{k \in [K]}\sum_{t \in [T]} c^-_{kt} \sigma_{kt} J\left(\frac{\epsilon^{\tinyU}_{kt} - \mu_{kt}}{\sqrt{3}\sigma_{kt}}\right) \\
= & g\left(\left\{\left(c^-_{kt}, \frac{\varepsilon^{\tinyU}_{kt} - \mu_{kt}}{\sqrt{3} \sigma_{kt}}\right)\right\}_{k \in [K], t \in [T]}\right).
\end{align*}
Similarly, $P^+(\varepsilon^{\tinyL}, \varepsilon^{\tinyU}) = g\Bigl(\left\{\left(c^+_{kt}, (\mu_{kt} - \varepsilon^{\tinyL}_{kt})/(\sqrt{3} \sigma_{kt})\right)\right\}_{k \in [K],}$ $_{t \in [T]}\Bigr)$.
%%%%%%%%%%%%%%%%%%%%%%%%%%%%%%%%%
\begin{figure}[h]
\centering
\includegraphics[width=0.7\linewidth]{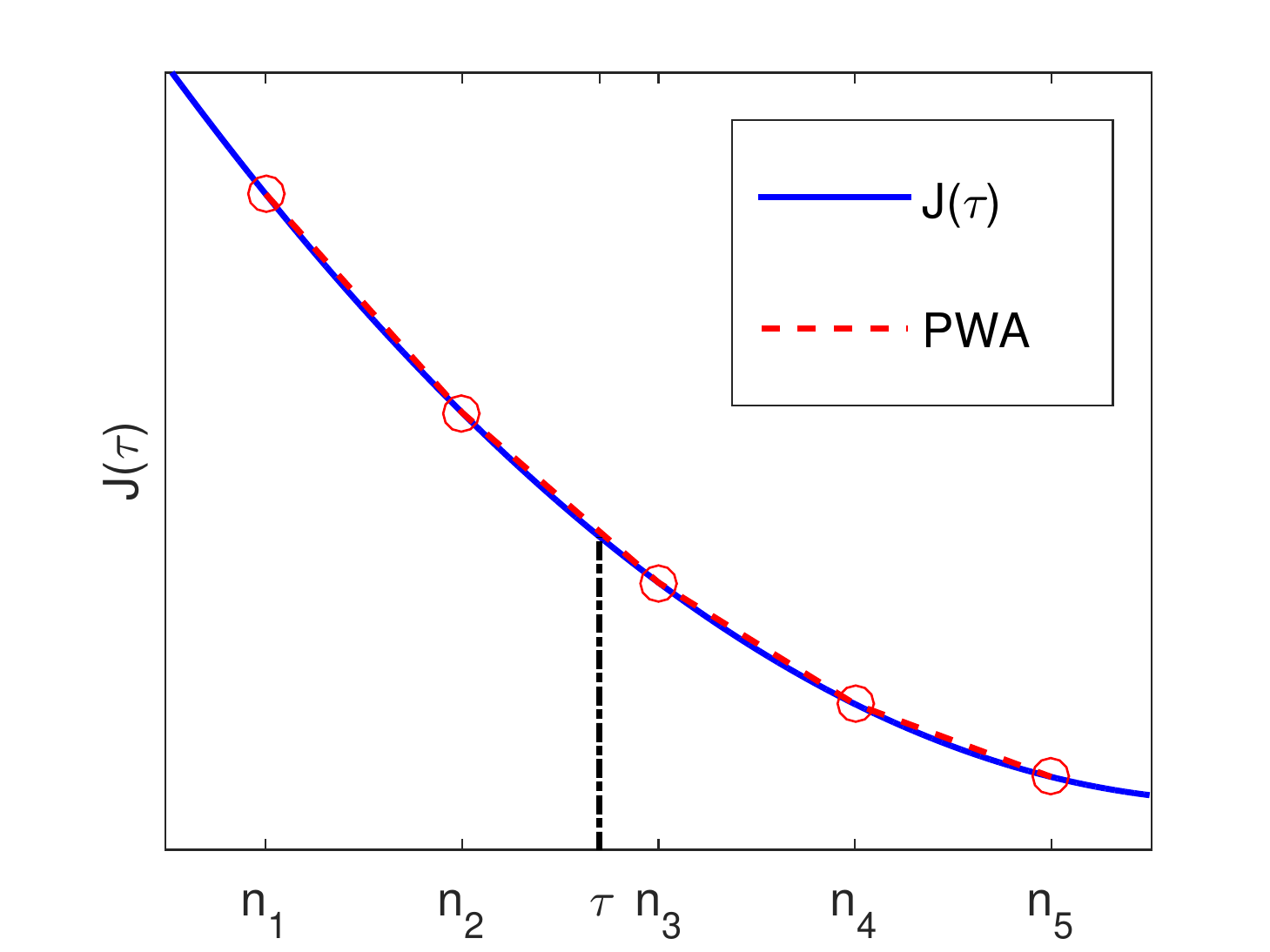}
\caption{An Illustration of $J(\tau)$ and its CPLA ($H = 4$)}
\label{pwa}
\end{figure}
%%%%%%%%%%%%%%%%%%%%%%%%%%%%%%%%%

Third, $J(\tau)$ is convex in $\tau$ because $\tau$ appears in the right-hand side of the formulation \eqref{apx-thm-2-note-12}--\eqref{apx-thm-2-note-13} with a minimization objective function that defines $J(\tau)$. It follows that $J(\tau)$ admits a conservative piecewise linear approximation (CPLA, see Fig.~\ref{pwa}). Specifically, suppose that $J(\tau)$ is defined on the interval $[\tau_{\tinyL}, \tau_{\tinyU}]$. Then, letting $n_h = \tau_{\tinyL} + (h-1)(\tau_{\tinyU} - \tau_{\tinyL})/H$ for $H \in \mathbb{N}_+$ and $h \in [H+1]$, we have $J(\tau) \leq \sum_{h=1}^{H+1} \lambda_h J(n_h)$, where
\begin{align*}
& \sum_{h=1}^{H+1} \lambda_h n_h = \tau, \\
& \sum_{h=1}^{H+1} \lambda_h = 1, \ \ \lambda_h \geq 0, \ \forall h \in [H+1].
\end{align*}
Hence, $J(\tau)$ can be conservatively approximated by $\sum_{h=1}^{H+1} \lambda_h J(n_h)$. Additionally, by construction, it is clear that $\lim_{H\rightarrow\infty} \sum_{h=1}^{H+1} \lambda_h J(n_h) = J(\tau)$. For all $k \in [K]$ and $t \in [T]$, we have $0 \leq \varepsilon^{\tinyU}_{kt} \leq w^{\text{max}}_{kt} - \hat{w}_{kt}$ and $w^{\text{min}}_{kt} - \hat{w}_{kt} \leq \varepsilon^{\tinyL}_{kt} \leq 0$ by constraint \eqref{MOc10}. It follows that $P^-(\varepsilon^{\tinyL}, \varepsilon^{\tinyU})$ can be conservatively approximated by $\sqrt{3} \sum_{k \in [K]} \sum_{t \in [T]} \sum_{h=1}^{H+1} c^-_{kt} \sigma_{kt} \lambda^-_{hkt} J(n^-_{hkt})$, where $n^-_{1kt} = -\mu_{kt}/(\sqrt{3}\sigma_{kt})$, $n^-_{(H+1)kt} = (w^{\text{max}}_{kt} - \hat{w}_{kt} - \mu_{kt})/(\sqrt{3}\sigma_{kt})$, $n^-_{hkt} = n^-_{1kt} + (h-1)(n^-_{(H+1)kt} - n^-_{1kt})/H$ for all $h \in [H+1]$, and
\begin{align*}
& \sum_{h=1}^{H+1} \lambda^-_{hkt} n^-_{hkt} = \frac{\varepsilon^{\tinyU}_{kt} - \mu_{kt}}{\sqrt{3} \sigma_{kt}}, \\
& \sum_{h=1}^{H+1} \lambda^-_{hkt} = 1, \ \ \lambda^-_{hkt} \geq 0, \ \forall h \in [H+1].
\end{align*}
Similarly, $P^+(\varepsilon^{\tinyL}, \varepsilon^{\tinyU})$ can be conservatively approximated by $\sqrt{3} \sum_{k \in [K]} \sum_{t \in [T]} \sum_{h=1}^{H+1} c^+_{kt} \sigma_{kt} \lambda^+_{hkt} J(n^+_{hkt})$, where $n^+_{1kt} = \mu_{kt}/(\sqrt{3}\sigma_{kt})$, $n^+_{(H+1)kt} = (\mu_{kt} + \hat{w}_{kt} - w^{\text{min}}_{kt})/(\sqrt{3}\sigma_{kt})$, $n^+_{hkt} = n^+_{1kt} + (h-1)(n^+_{(H+1)kt} - n^+_{1kt})/H$ for all $h \in [H+1]$, and
\begin{align*}
& \sum_{h=1}^{H+1} \lambda^+_{hkt} n^+_{hkt} = \frac{\mu_{kt} - \varepsilon^{\tinyL}_{kt}}{\sqrt{3} \sigma_{kt}}, \\
& \sum_{h=1}^{H+1} \lambda^+_{hkt} = 1, \ \ \lambda^+_{hkt} \geq 0, \ \forall h \in [H+1].
\end{align*}
\end{subequations}
Finally, we note that the values of $J(n^{\pm}_{hkt})$ can be efficiently obtained by solving the conic program \eqref{apx-thm-2-note-12}--\eqref{apx-thm-2-note-13} by setting $\tau := n^{\pm}_{hkt}$. Hence, we can incorporate the conservative approximations of $P^{\pm}(\varepsilon^{\tinyL}, \varepsilon^{\tinyU})$ into the DRCO model by using a set of linear constraints.

%\section*{Acknowledgment}

\bibliographystyle{IEEEtran}

\begin{thebibliography}{10}
\providecommand{\url}[1]{#1}
\csname url@samestyle\endcsname
\providecommand{\newblock}{\relax}
\providecommand{\bibinfo}[2]{#2}
\providecommand{\BIBentrySTDinterwordspacing}{\spaceskip=0pt\relax}
\providecommand{\BIBentryALTinterwordstretchfactor}{4}
\providecommand{\BIBentryALTinterwordspacing}{\spaceskip=\fontdimen2\font plus
\BIBentryALTinterwordstretchfactor\fontdimen3\font minus
  \fontdimen4\font\relax}
\providecommand{\BIBforeignlanguage}[2]{{%
\expandafter\ifx\csname l@#1\endcsname\relax
\typeout{** WARNING: IEEEtran.bst: No hyphenation pattern has been}%
\typeout{** loaded for the language `#1'. Using the pattern for}%
\typeout{** the default language instead.}%
\else
\language=\csname l@#1\endcsname
\fi
#2}}
\providecommand{\BIBdecl}{\relax}
\BIBdecl

\bibitem{wu2007stochastic}
L.~Wu, M.~Shahidehpour, and T.~Li, ``Stochastic security-constrained unit
  commitment,'' \emph{IEEE Transactions on Power Systems}, vol.~22, no.~2, pp.
  800--811, 2007.

\bibitem{papavasiliou2011stochastic}
A.~Papavasiliou, S.~S. Oren, and R.~P. O'Neill, ``Reserve requirements for wind
  power integration: A scenario-based stochastic programming framework,''
  \emph{IEEE Transactions on Power Systems}, vol.~26, no.~4, pp. 2197--2206,
  2011.

\bibitem{wang2012stochastic}
Q.~Wang, Y.~Guan, and J.~Wang, ``A chance-constrained two-stage stochastic
  program for unit commitment with uncertain wind power output,'' \emph{IEEE
  Transactions on Power Systems}, vol.~27, no.~1, pp. 206--215, 2012.

\bibitem{jiang2012robust}
R.~Jiang, J.~Wang, and Y.~Guan, ``Robust unit commitment with wind power and
  pumped storage hydro,'' \emph{IEEE Transactions on Power Systems}, vol.~27,
  no.~2, pp. 800--810, 2012.

\bibitem{bertsimas2013adaptive}
D.~Bertsimas, E.~Litvinov, X.~A. Sun, J.~Zhao, and T.~Zheng, ``Adaptive robust
  optimization for the security constrained unit commitment problem,''
  \emph{IEEE Transactions on Power Systems}, vol.~28, no.~1, pp. 52--63, 2013.

\bibitem{zhao2015dnelimit}
J.~Zhao, T.~Zheng, and E.~Litvinov, ``Variable resource dispatch through
  do-not-exceed limit,'' \emph{IEEE Transactions on Power Systems}, vol.~30,
  no.~2, pp. 820--828, 2015.

\bibitem{wei2015dispatchable}
W.~Wei, F.~Liu, and S.~Mei, ``Dispatchable region of the variable wind
  generation,'' \emph{IEEE Transactions on Power Systems}, vol.~30, no.~5, pp.
  2755--2765, 2015.

\bibitem{shao2017security}
C.~Shao, X.~Wang, M.~Shahidehpour, X.~Wang, and B.~Wang, ``Security-constrained
  unit commitment with flexible uncertainty set for variable wind power,''
  \emph{IEEE Transactions on Sustainable Energy}, vol.~8, no.~3, pp.
  1237--1246, 2017.

\bibitem{zymler2011distributionally}
S.~Zymler, D.~Kuhn, and B.~Rustem, ``Distributionally robust joint chance
  constraints with second-order moment information,'' \emph{Mathematical
  Programming}, vol. 137, no. 1-2, pp. 167--198, 2013.

\bibitem{wiesemann2014distributionally}
W.~Wiesemann, D.~Kuhn, and M.~Sim, ``Distributionally robust convex
  optimization,'' \emph{Operations Research}, vol.~62, no.~6, pp. 1358--1376,
  2014.

\bibitem{jiang2016data}
R.~Jiang and Y.~Guan, ``Data-driven chance constrained stochastic program,''
  \emph{Mathematical Programming}, vol. 158, no. 1-2, pp. 291--327, 2016.

\bibitem{zhao2016framework}
J.~Zhao, T.~Zheng, and E.~Litvinov, ``A unified framework for defining and
  measuring flexibility in power system,'' \emph{IEEE Transactions on Power
  Systems}, vol.~31, no.~1, pp. 339--347, 2016.

\bibitem{wei2015real}
W.~Wei, F.~Liu, and S.~Mei, ``Real-time dispatchability of bulk power systems
  with volatile renewable generations,'' \emph{IEEE Transactions on Sustainable
  Energy}, vol.~6, no.~3, pp. 738--747, 2015.

\bibitem{li2015robust}
Z.~Li, W.~Wu, B.~Zhang, and B.~Wang, ``Robust look-ahead power dispatch with
  adjustable conservativeness accommodating significant wind power
  integration,'' \emph{IEEE Transactions on Sustainable Energy}, vol.~6, no.~3,
  pp. 781--790, 2015.

\bibitem{li2015adjustable}
------, ``Adjustable robust real-time power dispatch with large-scale wind
  power integration,'' \emph{IEEE Transactions on Sustainable Energy}, vol.~6,
  no.~2, pp. 357--368, 2015.

\bibitem{wei2016dispatchability}
W.~Wei, J.~Wang, and S.~Mei, ``Dispatchability maximization for co-optimized
  energy and reserve dispatch with explicit reliability guarantee,'' \emph{IEEE
  Transactions on Power Systems}, vol.~31, no.~4, pp. 3276--3288, 2016.

\bibitem{shao2017power}
C.~Shao, X.~Wang, M.~Shahidehpour, X.~Wang, and B.~Wang, ``Power system
  economic dispatch considering steady-state secure region for wind power,''
  \emph{IEEE Transactions on Sustainable Energy}, vol.~8, no.~1, pp. 268--278,
  2017.

\bibitem{zeng2013solving}
B.~Zeng and L.~Zhao, ``Solving two-stage robust optimization problems using a
  column-and-constraint generation method,'' \emph{Operations Research
  Letters}, vol.~41, no.~5, pp. 457--461, 2013.

\bibitem{wang2017robustrcuc}
C.~Wang, F.~Liu, J.~Wang, F.~Qiu, W.~Wei, S.~Mei, and S.~Lei, ``Robust
  risk-constrained unit commitment with large-scale wind generation: An
  adjustable uncertainty set approach,'' \emph{IEEE Transactions on Power
  Systems}, vol.~32, no.~1, pp. 723--733, 2017.

\bibitem{wang2016risk}
C.~Wang, F.~Liu, J.~Wang, W.~Wei, and S.~Mei, ``Risk-based admissibility
  assessment of wind generation integrated into a bulk power system,''
  \emph{IEEE Transactions on Sustainable Energy}, vol.~7, no.~1, pp. 325--336,
  2016.

\bibitem{qiu2017data}
F.~Qiu, Z.~Li, and J.~Wang, ``A data-driven approach to improve wind
  dispatchability,'' \emph{IEEE Transactions on Power Systems}, vol.~32, no.~1,
  pp. 421--429, 2017.

\bibitem{li2016multi}
Z.~Li, F.~Qiu, and J.~Wang, ``Multi-period do-not-exceed limit for variable
  renewable generation dispatch considering discrete recourse controls,''
  \emph{arXiv preprint arXiv:1608.05273}, 2016.

\bibitem{korad2015zonal}
A.~S. Korad and K.~W. Hedman, ``Zonal do-not-exceed limits with robust
  corrective topology control,'' \emph{Electric Power Systems Research}, vol.
  129, pp. 235--242, 2015.

\bibitem{xiong2017distributionally}
P.~Xiong, P.~Jirutitijaroen, and C.~Singh, ``A distributionally robust
  optimization model for unit commitment considering uncertain wind power
  generation,'' \emph{IEEE Transactions on Power Systems}, vol.~32, no.~1, pp.
  39--49, 2017.

\bibitem{zhang2017distributionally}
Y.~Zhang, S.~Shen, and J.~L. Mathieu, ``Distributionally robust
  chance-constrained optimal power flow with uncertain renewables and uncertain
  reserves provided by loads,'' \emph{IEEE Transactions on Power Systems},
  vol.~32, no.~2, pp. 1378--1388, 2017.

\bibitem{xie2018distributionally}
W.~Xie and S.~Ahmed, ``Distributionally robust chance constrained optimal power
  flow with renewables: A conic reformulation,'' \emph{IEEE Transactions on
  Power Systems}, vol.~33, no.~2, pp. 1860--1867, 2018.

\bibitem{zhao2018distributionally}
C.~Zhao and R.~Jiang, ``Distributionally robust contingency-constrained unit
  commitment,'' \emph{IEEE Transactions on Power Systems}, vol.~33, no.~1, pp.
  94--102, 2018.

\bibitem{esfahani2017data}
P.~M. Esfahani and D.~Kuhn, ``Data-driven distributionally robust optimization
  using the wasserstein metric: Performance guarantees and tractable
  reformulations,'' \emph{Mathematical Programming}, Forthcoming, 2017.

\bibitem{wang2017risk}
C.~Wang, R.~Gao, F.~Qiu, J.~Wang, and L.~Xin, ``Risk-based distributionally
  robust optimal power flow with dynamic line rating,'' \emph{arXiv preprint
  arXiv:1712.08015}, 2017.

\bibitem{chen2018distributionally}
Y.~Chen, Q.~Guo, H.~Sun, Z.~Li, W.~Wu, and Z.~Li, ``A distributionally robust
  optimization model for unit commitment based on kullback-leibler
  divergence,'' \emph{IEEE Transactions on Power Systems}, Forthcoming, 2018.

\bibitem{li2016distributionally}
B.~Li, R.~Jiang, and J.~L. Mathieu, ``Distributionally robust risk-constrained
  optimal power flow using moment and unimodality information,'' in
  \emph{Decision and Control (CDC), 2016 IEEE 55th Conference on}.\hskip 1em
  plus 0.5em minus 0.4em\relax IEEE, 2016, pp. 2425--2430.

\bibitem{li2017ambiguous}
------, ``Ambiguous risk constraints with moment and unimodality information,''
  \emph{Mathematical Programming}, pp. 1--42, 2017.

\bibitem{fink2009wind}
S.~Fink, C.~Mudd, K.~Porter, and B.~Morgenstern, ``Wind energy curtailment case
  studies,'' National Renewable Energy Laboratory Report, Tech. Rep., 2009,
  available at http://www.nrel.gov/docs/fy10osti/46716.pdf.

\bibitem{doherty2005new}
R.~Doherty and M.~O'malley, ``A new approach to quantify reserve demand in
  systems with significant installed wind capacity,'' \emph{IEEE Transactions
  on Power Systems}, vol.~20, no.~2, pp. 587--595, 2005.

\bibitem{wang2008security}
J.~Wang, M.~Shahidehpour, and Z.~Li, ``Security-constrained unit commitment
  with volatile wind power generation,'' \emph{IEEE Transactions on Power
  Systems}, vol.~23, no.~3, pp. 1319--1327, 2008.

\bibitem{hodge2012wind}
B.-M. Hodge, D.~Lew, M.~Milligan, H.~Holttinen, S.~Sillanp{\"a}{\"a},
  E.~G{\'o}mez-L{\'a}zaro, R.~Scharff, L.~S{\"o}der, X.~G. Lars{\'e}n,
  G.~Giebel, D.~Flynn, and J.~Dobschinski, ``Wind power forecasting error
  distributions: {An} international comparison,'' in \emph{11th Annual
  International Workshop on Large-Scale Integration of Wind Power into Power
  Systems as well as on Transmission Networks for Offshore Wind Power Plants
  Conference}, 2012.

\bibitem{soyster1973convex}
A.~L. Soyster, ``Convex programming with set-inclusive constraints and
  applications to inexact linear programming,'' \emph{Operations Research},
  vol.~21, no.~5, pp. 1154--1157, 1973.

\bibitem{zimmerman2011matpower}
R.~D. Zimmerman, C.~E. Murillo-S{\'a}nchez, and R.~J. Thomas, ``Matpower:
  Steady-state operations, planning, and analysis tools for power systems
  research and education,'' \emph{IEEE Transactions on power systems}, vol.~26,
  no.~1, pp. 12--19, 2011.

\bibitem{nreleasternwind}
``{NREL Eastern Wind Dataset},''
  \url{https://www.nrel.gov/grid/eastern-wind-data.html}, accessed: 2018-07-30.

\bibitem{xie2017optimized}
W.~Xie, S.~Ahmed, and R.~Jiang, ``Optimized bonferroni approximations of
  distributionally robust joint chance constraints,'' \emph{Available at
  Optimization Online}, 2017.

\bibitem{gauss1823theoria}
C.-F. Gauss, \emph{Theoria combinationis observationum erroribus minimis
  obnoxiae}.\hskip 1em plus 0.5em minus 0.4em\relax Henricus Dieterich, 1823,
  vol.~1.

\bibitem{ben2001lectures}
A.~Ben-Tal and A.~Nemirovski, \emph{Lectures on modern convex optimization:
  analysis, algorithms, and engineering applications}.\hskip 1em plus 0.5em
  minus 0.4em\relax SIAM, 2001, vol.~2.

\bibitem{dharmadhikari1988unimodality}
S.~Dharmadhikari and K.~Joag-Dev, \emph{Unimodality, Convexity, and
  Applications}.\hskip 1em plus 0.5em minus 0.4em\relax Elsevier, 1988.

\bibitem{bertsimas2005optimal}
D.~Bertsimas and I.~Popescu, ``Optimal inequalities in probability theory: A
  convex optimization approach,'' \emph{SIAM Journal on Optimization}, vol.~15,
  no.~3, pp. 780--804, 2005.

\end{thebibliography}
\end{document}